
\def\LaTeX{%
  \let\Begin\begin
  \let\End\end
  \let\salta\relax
  \let\finqui\relax
  \let\futuro\relax}

\def\UK{\def\our{our}\let\sz s}
\def\USA{\def\our{or}\let\sz z}

\UK



\LaTeX

\USA


\salta

\documentclass[twoside,12pt]{article}
\setlength{\textheight}{24cm}
\setlength{\textwidth}{16cm}
\setlength{\oddsidemargin}{2mm}
\setlength{\evensidemargin}{2mm}
\setlength{\topmargin}{-15mm}
\parskip2mm


\usepackage[usenames,dvipsnames]{color}
\usepackage{amsmath}
\usepackage{amsthm}
\usepackage{amssymb}
\usepackage[mathcal]{euscript}

%
%


\definecolor{viola}{rgb}{0.3,0,0.7}
\definecolor{ciclamino}{rgb}{0.5,0,0.5}

\def\pier #1{{\color{red}#1}}
\def\juerg #1{{\color{Green}#1}}
\def\gianni #1{{\color{blue}#1}}
\def\pier #1{#1}
\def\juerg #1{#1}
\def\gianni #1{#1}



\bibliographystyle{plain}


%

\finqui

\def\Beq{\Begin{equation}}
\def\Eeq{\End{equation}}
\def\Bsist{\Begin{eqnarray}}
\def\Esist{\End{eqnarray}}

\def\Bthm{\Begin{theorem}}
\def\Ethm{\End{theorem}}

\def\Bprop{\Begin{proposition}}
\def\Eprop{\End{proposition}}
\def\Bcor{\Begin{corollary}}
\def\Ecor{\End{corollary}}
\def\Brem{\Begin{remark}\rm}
\def\Erem{\End{remark}}

\def\Bdim{\Begin{proof}}
\def\Edim{\End{proof}}
\def\Bcenter{\Begin{center}}
\def\Ecenter{\End{center}}
\let\non\nonumber




\def\step #1 \par{\medskip\noindent{\bf #1.}\quad}


\def\Lip{Lip\-schitz}

\def\Frechet{Fr\'echet}
\def\aand{\quad\hbox{and}\quad}

\def\lhs{left-hand side}
\def\rhs{right-hand side}
\def\sfw{straightforward}


\def\generaliz{generali\sz}

\def\minimiz{minimi\sz}

\def\recogniz{recogni\sz}


\def\multibold #1{\def\arg{#1}%
  \ifx\arg\pto \let\next\relax
  \else
  \def\next{\expandafter
    \def\csname #1#1#1\endcsname{{\bf #1}}%
    \multibold}%
  \fi \next}

\def\pto{.}

\def\multical #1{\def\arg{#1}%
  \ifx\arg\pto \let\next\relax
  \else
  \def\next{\expandafter
    \def\csname cal#1\endcsname{{\cal #1}}%
    \multical}%
  \fi \next}


\def\multimathop #1 {\def\arg{#1}%
  \ifx\arg\pto \let\next\relax
  \else
  \def\next{\expandafter
    \def\csname #1\endcsname{\mathop{\rm #1}\nolimits}%
    \multimathop}%
  \fi \next}

\multibold
qwertyuiopasdfghjklzxcvbnmQWERTYUIOPASDFGHJKLZXCVBNM.

\multical
QWERTYUIOPASDFGHJKLZXCVBNM.

\multimathop
ad dist div dom meas sign supp .


\def\accorpa #1#2{\eqref{#1}--\eqref{#2}}
\def\Accorpa #1#2 #3 {\gdef #1{\eqref{#2}--\eqref{#3}}%
  \wlog{}\wlog{\string #1 -> #2 - #3}\wlog{}}


\def\separa{\noalign{\allowbreak}}

\def\neto{\mathrel{{\scriptscriptstyle\nearrow}}}
\def\seto{\mathrel{{\scriptscriptstyle\searrow}}}

\def\graffe #1{\mathopen\{#1\mathclose\}}

\def\<#1>{\mathopen\langle #1\mathclose\rangle}
\def\norma #1{\mathopen \| #1\mathclose \|}

\def\[#1]{\mathopen\langle\!\langle #1\mathclose\rangle\!\rangle}

\def\iot {\int_0^t}
\def\ioT {\int_0^T}
\def\intQt{\int_{Q_t}}
\def\intQ{\int_Q}
\def\iO{\int_\Omega}
\def\iG{\int_\Gamma}
\def\intS{\int_\Sigma}
\def\intSt{\int_{\Sigma_t}}

\def\dt{\partial_t}
\def\dn{\partial_n}

\def\cpto{\,\cdot\,}

\def\checkmmode #1{\relax\ifmmode\hbox{#1}\else{#1}\fi}

\def\aeQ{\checkmmode{a.e.\ in~$Q$}}

\def\aeS{\checkmmode{a.e.\ on~$\Sigma$}}
\def\aet{\checkmmode{a.e.\ in~$(0,T)$}}

\def\aat{\checkmmode{for a.a.~$t\in(0,T)$}}


\def\erre{{\mathbb{R}}}




\def\genspazio #1#2#3#4#5{#1^{#2}(#5,#4;#3)}
\def\spazio #1#2#3{\genspazio {#1}{#2}{#3}T0}

\def\L {\spazio L}
\def\H {\spazio H}

\def\C #1#2{C^{#1}([0,T];#2)}


\def\Lx #1{L^{#1}(\Omega)}
\def\Hx #1{H^{#1}(\Omega)}

\def\LxG #1{L^{#1}(\Gamma)}
\def\HxG #1{H^{#1}(\Gamma)}

\def\LQ #1{L^{#1}(Q)}
\def\LS #1{L^{#1}(\Sigma)}

\def\Ldue{\Lx 2}

\def\Huno{\Hx 1}
\def\Hdue{\Hx 2}

\def\HunoG{\HxG 1}
\def\HdueG{\HxG 2}

\def\LdueG{\LxG 2}


\def\LQ #1{L^{#1}(Q)}


\let\theta\vartheta

\let\phi\varphi
\let\lam\lambda

\let\TeXchi\chi                         
\newbox\chibox
\setbox0 \hbox{\mathsurround0pt $\TeXchi$}
\setbox\chibox \hbox{\raise\dp0 \box 0 }
\def\chi{\copy\chibox}



\def\suG{{\vrule height 5pt depth 4pt\,}_\Gamma}

\def\fG{f_\Gamma}
\def\yG{y_\Gamma}
\def\uG{u_\Gamma}
\def\vG{v_\Gamma}
\def\hG{h_\Gamma}
\def\xiG{\xi_\Gamma}
\def\qG{q_\Gamma}
\def\uGmin{u_{\Gamma,{\rm min}}}
\def\uGmax{u_{\Gamma,{\rm max}}}

\def\yt{y^\tau}
\def\ytG{y^\tau_\Gamma}
\def\wt{w^\tau}

\def\yz{y_0}

\def\Mz{M_0}
\def\ustar{u_*}
\def\vstar{v_*}

\def\bQ{b_Q}
\def\bS{b_\Sigma}

\def\bz{b_0}

\def\rmin{r_-}
\def\rmax{r_+}

\def\zQ{z_Q}
\def\zS{z_\Sigma}

\def\phQ{\phi_Q}
\def\phS{\phi_\Sigma}
\def\phtQ{\phi_Q^\tau}
\def\phtS{\phi_\Sigma^\tau}

\def\Uad{\calU_{\rm ad}}
\def\uopt{\overline u_\Gamma}
\def\yopt{\overline y}
\def\yoptG{\overline y_\Gamma}
\def\wopt{\overline w}

\def\utopt{\overline u^{\,\pier{\tau}}_\Gamma}
\def\ytopt{\overline y^{\,\pier{\tau}}}
\def\ytoptG{\overline y^{\,\pier{\tau}}_\Gamma}
\def\wtopt{\overline w^{\,\pier{\tau}}}

\def\utmod{\tilde u^\tau_\Gamma}
\def\ytmod{\tilde y^\tau}
\def\ytmodG{\tilde y^\tau_\Gamma}
\def\wtmod{\tilde w^\tau}

\def\utG{u^\tau_\Gamma}

\def\pt{p^\tau}
\def\qt{q^\tau}
\def\qtG{q^\tau_\Gamma}

\def\unG{u_{\Gamma,n}}
\def\yn{y_n}
\def\ynG{y_{\Gamma,n}}
\def\wn{w_n}

\def\uO{u^\Omega}
\def\vO{v^\Omega}

\def\VO{\calV_\Omega}
\def\HO{\calH_\Omega}

\def\calWz{\calW_0}
\def\calWp{\calW_0^*}

\def\VG{V_\Gamma}
\def\HG{H_\Gamma}
\def\Vp{V^*}
\def\VGp{\VG^*}

\def\normaV #1{\norma{#1}_V}
\def\normaH #1{\norma{#1}_H}

\def\normaVp #1{\norma{#1}_*}

\def\nablaG{\nabla_{\!\Gamma}}
\def\DeltaG{\Delta_\Gamma}

\let\hat\widehat

\def\Pi{\hat\pi}

\let\lam\lambda
\def\lamG{\lam_\Gamma}
\let\Lam\Lambda
\def\LamG{\Lam_\Gamma}

\def\lamt{\lambda^\tau}
\def\lamtG{\lambda^\tau_\Gamma}


\Begin{document}


\title{A boundary control problem\\[0.3cm] 
  for the pure Cahn--Hilliard equation\\[0.3cm]
  with dynamic boundary conditions\footnote{{\bf Acknowledgment.}\quad\rm
PC and GG gratefully acknowledge some financial support from the MIUR-PRIN Grant 2010A2TFX2 ``Calculus of Variations'' and the GNAMPA (Gruppo Nazionale per l'Analisi Matematica, la Probabilit\`a e le loro Applicazioni) of INdAM (Istituto Nazionale di Alta Matematica).}}

\author{}
\date{}
\maketitle
\Bcenter
\vskip-2cm
{\large\sc Pierluigi Colli$^{(1)}$}\\
{\normalsize e-mail: {\tt pierluigi.colli@unipv.it}}\\[.25cm]
{\large\sc Gianni Gilardi$^{(1)}$}\\
{\normalsize e-mail: {\tt gianni.gilardi@unipv.it}}\\[.25cm]
{\large\sc J\"urgen Sprekels$^{(2)}$}\\
{\normalsize e-mail: {\tt sprekels@wias-berlin.de}}\\[.45cm]
$^{(1)}$
{\small Dipartimento di Matematica ``F. Casorati'', Universit\`a di Pavia}\\
{\small via Ferrata 1, 27100 Pavia, Italy}\\[.2cm]
$^{(2)}$
{\small Weierstrass Institute for Applied Analysis and Stochastics}\\
{\small Mohrenstrasse 39, 10117 Berlin, Germany}\\[2mm]
{\small and}\\[2mm]
{\small Department of Mathematics}\\
{\small Humboldt-Universit\"at zu Berlin}\\
{\small \pier{Unter den Linden 6, 10099 Berlin, Germany}}\\
[1cm]
\Ecenter

\Begin{abstract}
A boundary control problem for the pure Cahn--Hilliard equations 
with possibly singular potentials and dynamic boundary conditions
is studied and \pier{first-order} necessary conditions for optimality are proved.
\vskip3mm

\noindent {\bf Key words:}
Cahn--Hilliard equation, dynamic boundary conditions, phase separation,
singular potentials, optimal control, optimality conditions.
\vskip3mm
\noindent {\bf AMS (MOS) Subject Classification:} 35K55 (35K50, 82C26)
\End{abstract}

\salta

\pagestyle{myheadings}
\newcommand\testopari{\sc Colli \ --- \ Gilardi \ --- \ Sprekels}
\newcommand\testodispari{\sc Boundary control problem for the pure Cahn--Hilliard equation}
\markboth{\testodispari}{\testopari}

\finqui


\section{Introduction}
\label{Intro}
\setcounter{equation}{0}

The simplest form of the Cahn--Hilliard equation (see \cite{CahH, EllSt, EllSh}) reads as follows
\Beq
  \dt y - \Delta w = 0
  \aand
  w = -\Delta y + f'(y) 
  \quad \hbox{in $\Omega\times (0,T) $},
  \label{Icahil}
\Eeq
where $\Omega$ is the domain where the evolution takes place,
and $y$~and $w$ denote the order parameter and the chemical potential, respectively. 
Moreover, $f'$ represents the derivative of a double well potential~$f$,
and typical and important examples are the following
\Bsist
  & f_{reg}(r) = \frac14(r^2-1)^2 \,,
  \quad r \in \erre 
  \label{regpot}
  \\
  & f_{log}(r) = ((1+r)\ln (1+r)+(1-r)\ln (1-r)) - c r^2 \,,
  \quad r \in (-1,1),
  \label{logpot}
\Esist
where $c>0$ in the latter is large enough in order \pier{that $f_{log}$ 
be} nonconvex. 
The potentials \eqref{regpot} and \eqref{logpot}
are usually called the classical regular potential
and the logarithmic double-well potential, respectively.

The present paper is devoted to the study of the control problem described below
for the initial--boundary value problem
obtained by complementing \eqref{Icahil} with an initial condition like $y(0)=\yz$ 
and the following boundary conditions
\Beq
  \dn w = 0 
  \aand
  \dn y + \dt\yG - \DeltaG\yG + \fG'(\yG) = \uG
  \quad \hbox{on $\Gamma\times (0,T)$}
  \label{Ibc}
\Eeq
where $\Gamma$ is the boundary of~$\Omega$.
The former is very common in the literature and preserves mass conservation,
i.e., it implies that the space integral of $y$ is constant in time.
The latter is an evolution equation
for the trace $\yG$ of the order parameter on the boundary,
and the normal derivative $\dn y$ and $\uG$ act as forcing terms.
\pier{This condition enters the class of the so-called \gianni{dynamic boundary conditions
that have been widely used} in the literature in the last
twenty years, say: in particular, the study of dynamic boundary conditions 
with Cahn--Hilliard type equations has been taken up by some 
authors (let us quote \cite{CMZ, CGS, GiMiSchi, PRZ, RZ, WZ} and also refer to 
the recent contribution \cite{CF} in which also a forced mass constraint 
on the boundary is considered).}

\pier{The dynamic boundary condition in \eqref{Ibc} contains 
the Laplace-Beltrami operator $\DeltaG$ and a
nonlinearity $\fG'$ which is analogous to $f'$ 
but is} now acting on the boundary values~$\uG$.
Even though some of our \juerg{results} hold under weaker hypotheses,
we assume from the very beginning that $f'$ and $\fG'$ have the same domain~$\calD$.
The main assumption we make is a compatibility condition between 
\juerg{these} nonlinearities.
Namely, we suppose that $\fG'$ dominates $f'$ in the following sense:
\Beq
  |f'(r)| \leq \eta \, |\fG'(r)| + C
  \label{Icompatib}
\Eeq
for some positive constants $\eta$ and~$C$ 
and for every $r\in\calD$.
\juerg{This} condition, earlier introduced in~\cite{CaCo} 
in relation with the Allen--Cahn equation with dynamic boundary conditions \pier{(see also \cite{CS}),} is then used in~\cite{CGS} (as~well as in \cite{CFGS1} and~\cite{CGSvisc})
to deal with the Cahn--Hilliard system.
This complements~\cite{GiMiSchi}, 
where \juerg{some} kind of an opposite inequality is assumed.

As just said, this paper deals with a control problem for the state system
described above, the control being the source term $\uG$ that appears 
in the dynamic boundary condition~\eqref{Ibc}.
Namely, the problem we want to address consists in \minimiz ing
a proper cost functional depending on both the control $\uG$
and the associate state $(y,\yG)$.
Among several possibilities, we choose the cost functional 
\Beq
  \calJ(y,\yG,\uG)
  := \frac\bQ 2 \, \norma{y-\zQ}_{L^2(Q)}^2
  + \frac\bS 2 \, \norma{\yG-\zS}_{L^2(\Sigma)}^2
  + \frac\bz 2 \, \norma\uG_{L^2(\Sigma)}^2\,,
  \label{Icost}
\Eeq
where the functions $z_Q, z_\Sigma$ and the nonnegative constants $\bQ,\bS,\bz$ are given. 
The control problem then consists in \minimiz ing \eqref{Icost} 
subject to the state system and to the constraint $\uG\in\Uad$, 
where the control box $\Uad$ is given~by
\Bsist
  & \Uad := 
  & \bigl\{ \uG\in\H1\HG\cap\LS\infty:
  \non
  \\
  && \ \uGmin\leq\uG\leq\uGmax\ \aeS,\ \norma{\dt\uG}_{\LS2}\leq\Mz
  \bigr\} 
  \label{Iuad}
\Esist
for some given functions $\uGmin,\uGmax\in\LS\infty$
and some prescribed positive constant~$\Mz$. \pier{Of course, the control box  $\Uad$ must be nonempty and this is guaranteed if, for instance, at least one of  $\uGmin$ or $\uGmax $ is in $\H1\HG$} \juerg{and its time derivative satisfies the above $L^2(\Sigma)$ bound.}

This paper \pier{is a follow-up of} the recent 
contributions \cite{CGS} and \cite{CGSvisc} already mentioned.
They deal with a similar system and a similar control problem.
\pier{The paper} \cite{CGS}~contains a number of results on the state system
obtained by considering
\Beq
  w = \tau \, \dt y - \Delta y + f'(y)
  \label{Ivisc}
\Eeq
in place of the second \pier{condition in}~\eqref{Icahil}.
In~\eqref{Ivisc}, $\tau$~is a nonnegative parameter
and the case $\tau>0$ coupled with the first \pier{equation in} 
\eqref{Icahil} yields the well-known viscous Cahn--Hilliard equation
(in~contrast, we term \eqref{Icahil} the pure Cahn--Hilliard system).
More precisely, existence, uniqueness and regularity results are proved in \cite{CGS}
for general potentials that include \accorpa{regpot}{logpot},
and are valid for both the viscous and pure cases, i.e., by assuming just $\tau\geq0$.
Moreover, if $\tau>0$, further regularity and properties of the solution are ensured.
\juerg{These} results are then used in \cite{CGSvisc}, where the boundary control problem
associated to a cost functional that \generaliz es \eqref{Icost} 
is~addressed
and both the existence of an optimal control and \pier{first-order} necessary \pier{conditions} for optimality \pier{are proved and expressed}
in terms of the solution of a proper adjoint problem.

\pier{In fact, recently Cahn--Hilliard systems have been rather investigated from the viewpoint of optimal control. In this connection, we refer to \cite{HW1, WaNa, ZW} and to \cite{ZL1, ZL2} which deal with the convective Cahn--Hilliard equation; 
the case with a nonlocal potential \pier{is} studied in \cite{RoSp}. There also exist contributions addressing some discretized versions of general Cahn--Hilliard systems, cf. \cite{HW2, Wang}. However, about the optimal control of viscous or non-viscous Cahn--Hilliard systems with dynamic boundary conditions of the form (\ref{Ibc}), we only know of the papers 
\cite{CGSvisc} and \cite{CFGS1} dealing with the viscous case; to the best of our knowledge, the present contribution is the first paper treating 
the optimal control of the pure Cahn--Hilliard system with dynamic boundary conditions.}

The technique used in \pier{our approach} essentially consists in
starting from the known results for $\tau>0$
and then letting the parameter $\tau$ tend to zero.
In doing that, we use some of the ideas of \cite{CFS} and \cite{CFGS1},
which deal with the Allen--Cahn and the viscous Cahn--Hilliard equations, respectively,
and address similar control problems related to the nondifferentiable
double obstacle potential by seeing it as a limit of logarithmic double-well potentials.

The paper is organized as follows.
In the next section, we list our assumptions, state the problem in a precise form
and present our results.
The corresponding proofs are given in the last section.


\section{Statement of the problem and results}
\label{STATEMENT}
\setcounter{equation}{0}

In this section, we describe the problem under study and give an outline of our results.
As in the Introduction,
$\Omega$~is the body where the evolution takes place.
We assume $\Omega\subset\erre^3$
to~be open, bounded, connected, and smooth,
and we write $|\Omega|$ for its Lebesgue measure.
Moreover, $\Gamma$, $\dn$, $\nablaG$ and $\DeltaG$ stand for
the boundary of~$\Omega$, the outward normal derivative, the surface gradient 
and the Laplace--Beltrami operator, respectively.
Given a finite final time~$T>0$,
we set for convenience
\Bsist
  && Q_t := \Omega \times (0,t)
  \aand
  \Sigma_t := \Gamma \times (0,t)
  \quad \hbox{for every $t\in(0,T]$}
  \label{defQtSt}
  \\
  && Q := Q_T \,,
  \aand
  \Sigma := \Sigma_T \,.
  \label{defQS}
\Esist
\Accorpa\defQeS defQtSt defQS
Now, we specify the assumptions on the structure of our system.
Even though some of the results we quote hold under rather mild hypotheses,
we give a list of assumptions that implies the whole set of conditions required in~\cite{CGS}.
We assume~that
\Bsist
  & -\infty \leq \rmin < 0 < \rmax \leq +\infty
  \label{hpDf}\\[1mm]
    & \hbox{$f,\,\fG:(\rmin,\rmax)\to[0,+\infty)$
    are $C^3$ functions}
  \label{hppot}
  \\[1mm]
  & f(0) = \fG(0) = 0
  \aand 
  \hbox{$f''$ and $\fG''$ are bounded from below}
  \label{hpfseconda}
  \\[1mm]
  \separa
  & |f'(r)| \leq \eta \,|\fG'(r)| + C
  \quad \hbox{for some $\eta,\, C>0$ and every $r\in(\rmin,\rmax)$}
  \label{hpcompatib}
  \\[1mm]
  \separa
  & \lim\limits_{r\seto\rmin} f'(r)
  = \lim\limits_{r\seto\rmin} \fG'(r) 
  = -\infty 
  \aand
  \lim\limits_{r\neto\rmax} f'(r)
  = \lim\limits_{r\neto\rmax} \fG'(r) 
  = +\infty \, .
  \label{fmaxmon}
\Esist
\Accorpa\HPstruttura hpDf fmaxmon
We note that \HPstruttura\ imply the possibility
of splitting $f'$ as $f'=\beta+\pi$,
where $\beta$ is a monotone function that diverges at~$r_\pm$
and $\pi$ is a perturbation with a bounded derivative.
Moreover, the same is true for~$\fG$,
so that the assumptions of \cite{CGS} are satisfied.
Furthermore, the choices $f=f_{reg}$ and $f=f_{log}$
corresponding to \eqref{regpot} and \eqref{logpot} are allowed.

Next, in order to simplify notations, we~set
\Bsist
  && V := \Huno, \quad
  H := \Ldue, \quad
  \HG := \LdueG 
  \aand
  \VG := \HunoG 
  \qquad
  \label{defVH}
  \\
  && \calV := \graffe{(v,\vG)\in V\times\VG:\ \vG=v\suG}
  \aand
  \calH := H \times \HG   
  \label{defcalVH}
\Esist
\Accorpa\Defspazi defVH defcalVH
and endow these spaces with their natural norms.
If $X$ is any Banach space, \juerg{then} $\norma\cpto_X$~and $X^*$ denote 
its norm and its dual space, respectively.
Furthermore, the symbol $\<\cpto,\cpto>$ usually stands
for the duality pairing between~$\Vp$ and $V$ itself
and the similar notation $\<\cpto,\cpto>_\Gamma$ refers to the spaces $\VGp$ and~$\VG$.
In the following, it is understood that $H$ is identified with $H^*$ and thus embedded in $\Vp$
in the usual way, i.e., such that we have
$\<u,v>=(u,v)$ with the inner product $(\,\cdot\,,\,\cdot)$ of~$H$,
for every $u\in H$ and $v\in V$.
Thus, we introduce the Hilbert triplet $(V,H,\Vp)$
and analogously behave with the boundary spaces $\VG$, $\HG$ and~$\VGp$.
Finally, if $u\in\Vp$ and $\underline u\in\L1\Vp$,
we define their \generaliz ed mean values 
$\uO\in\erre$ and $\underline u^\Omega\in L^1(0,T)$ by setting
\Beq
  \uO := \frac 1 {|\Omega|} \, \< u , 1 >
  \aand
  \underline u^\Omega(t) := \bigl( \underline u(t) \bigr)^\Omega
  \quad \aat .
  \label{media}
\Eeq
Clearly, \juerg{the relations in} \eqref{media} give the usual mean values when applied to elements
of~$H$ or $\L1H$.

At this point, we can describe the state problem.
For the data, we assume~that
\Bsist
  && \yz \in \Hx 2 
  \aand
  \yz\suG \in \HxG 2
  \label{hpyz}
  \\
  && \rmin < \yz(x) < \rmax 
  \quad \hbox{for every $x\in\overline\Omega$}
  \label{hpyzbis}
  \\
  && \uG \in \H1\HG \,.
  \label{hpuG}
\Esist
\Accorpa\HPdati hpyz hpuG
We look for a triplet $(y,\yG,w)$ satisfying
\Bsist
  & y \in \H1\Vp \cap \L\infty V \cap \L2\Hdue
  \label{regy}
  \\
  & \yG \in \H1\HG \cap \L\infty\VG \cap \L2\HdueG
  \label{regyG}
  \\
  & \yG(t) = y(t)\suG
  \quad \aat
  \label{tracciay}
  \\
  & w \in \L2V\,,
  \label{regw}
\Esist
\Accorpa\Regsoluz regy regw
as well as, for almost every $t\in (0,T)$, the variational equations
\Bsist
  && \< \dt y(t) \, v >
  + \iO \nabla w(t) \cdot \nabla v = 0
  \qquad \hbox{for every $v\in V$}
  \label{prima}
  \\
  \noalign{\smallskip}
  && \iO w(t) \, v
  = \iO \nabla y(t) \cdot \nabla v
  + \iG \dt\yG(t) \, \vG
  + \iG \nablaG\yG(t) \cdot \nablaG\vG
  \qquad
  \non
  \\
  && \quad {}
  + \iO f'(y(t)) \, v
  + \iG \bigl( \fG'(\yG(t)) - \uG(t) \bigr) \, \vG
  \qquad \hbox{for every $(v,\vG)\in\calV$}
  \qquad
  \label{seconda}
  \\
  && y(0) = \yz \,. \vphantom\sum
  \label{cauchy}
\Esist
\Accorpa\State prima cauchy
\Accorpa\Pbl regy cauchy
Thus, we require that the state variables satisfy the variational counterpart 
of the problem described in the Introduction in a strong form.
We note that an equivalent formulation of \accorpa{prima}{seconda}
is given~by
\Bsist
  && \iot \< \dt y(t) \, v(t) > \, dt
  + \intQ \nabla w \cdot \nabla v = 0
  \label{intprima}
  \\
  \noalign{\smallskip}
  && \intQ wv
  = \intQ \nabla y \cdot \nabla v
  + \intS \dt\yG \, \vG
  + \intS \nablaG\yG \cdot \nablaG\vG
  \qquad
  \non
  \\
  && \quad {}
  + \intQ f'(y) \, v
  + \intS \bigl( \fG'(\yG) - \uG \bigr) \, \vG
  \qquad
  \label{intseconda}
\Esist
\Accorpa\IntPbl intprima intseconda
for every $v\in\L2 V$ and every $(v,\vG)\in\L2\calV$, respectively.

Besides, we consider the analogous state system with viscosity.
Namely, for $\tau>0$ we replace \eqref{seconda}~by
\Bsist
  && \iO w(t) \, v
  = \tau \iO \dt y(t) \, v
  + \iO \nabla y(t) \cdot \nabla v
  + \iG \dt\yG(t) \, \vG
  + \iG \nablaG\yG(t) \cdot \nablaG\vG
  \qquad
  \non
  \\
  && \quad {}
  + \iO f'(y(t)) \, v
  + \iG \bigl( \fG'(\yG(t)) - \uG(t) \bigr) \, \vG
  \qquad \hbox{for every $(v,\vG)\in\calV$}
  \qquad
  \label{secondavisc}
\Esist
\def\Pbltau{\eqref{prima}, \eqref{cauchy} and~\eqref{secondavisc}}%
in the above system.
We notice that a variational equation equivalent to \eqref{secondavisc} 
is given by the analogue of \eqref{intseconda}, i.e.,
\Bsist
  && \intQ wv
  = \tau \intQ \dt y \, v
  + \intQ \nabla y \cdot \nabla v
  + \intS \dt\yG \, \vG
  + \intS \nablaG\yG \cdot \nablaG\vG
  \qquad
  \non
  \\
  && \quad {}
  + \intQ f'(y) \, v
  + \intS \bigl( \fG'(\yG) - \uG \bigr) \, \vG
  \quad
  \hbox{for every $(v,\vG)\in\L2\calV$}.
  \qquad
  \label{intsecondavisc}
\Esist

As far as existence, uniqueness, regularity and continuous dependence are concerned,
we directly refer to~\cite{CGS}.
From \cite[Thms.~2.2 and~2.3]{CGS}
(where $\calV$ has a slightly different meaning with respect to the present paper),
we~have the following results:

\Bthm
\label{daCGSpure}
Assume \HPstruttura\ and \HPdati.
Then, there exists a unique triplet $(y,\yG,w)$ satisfying
\Regsoluz\ and solving \State.
\Ethm

\Bthm
\label{daCGS}
Assume \HPstruttura\ and \HPdati.
Then, for every $\tau>0$,
there exists a unique triplet $(\yt,\ytG,\wt)$ satisfying
\Regsoluz\ and solving \Pbltau.
Moreover, \juerg{this} solution satisfies
$\dt\yt\in\L2H$
and the estimate
\Bsist
  && \norma\yt_{\H1\Vp \cap \L\infty V \cap \L2\Hdue}
  \non
  \\[1mm]
  && \quad {}
  + \norma\ytG_{\H1\HG \cap \L\infty\VG \cap \L2\HdueG}
  \non
  \\[1mm]
  && \quad {}
  + \norma\wt_{\pier{\L2V}}
  + \norma{f'(\yt)}_{\L2H}
  + \norma{\fG'(\ytG)}_{\L2\HG}
  \non
  \\[1mm]
  && \quad {}
  + \tau^{1/2} \norma{\dt\yt}_{\L2H}
  \leq \, C_0 
  \label{stab}
\Esist
holds true for some constant $C_0>0$  that depends only on
$\Omega$, $T$, the shape of the nonlinearities $f$ and~$\fG$,
\pier{and the norms~$\bigl\Vert (\yz,{\yz}_{|_\Gamma})\bigr\Vert_{\calV}$,
$\Vert f' (\yz)\Vert_{L^1(\Omega)}$,
$\bigl\Vert \fG' \bigl({\yz}_{|_\Gamma}\bigr) \bigr\Vert_{L^1(\Gamma)}$,
and $\norma{\uG}_{\L2\HG}.$}
\Ethm

In fact, \pier{if the data are more regular, in particular\juerg{, if } 
$\uG \in \H1\HG \cap\LS\infty $,} 
\juerg{then} the solution $(\yt,\ytG,\wt)$ is even smoother
(see \cite[Thms.~2.4 and~2.6]{CGS}) \pier{and, specifically,}
it satisfies
\Beq
  \rmin^\tau \leq \yt \leq \rmax^\tau
  \quad \aeQ
  \label{faraway}
\Eeq
for some constants $\rmin^\tau,\rmax^\tau\in(\rmin,\rmax)$
that depend on~$\tau$, in addition.
It \juerg{follows} that the functions
$f''(\yt)$ and $\fG''(\ytG)$ 
(which will appear as coefficients in a linear system later~on) 
are~bounded.
We also notice that the stability estimate \eqref{stab} is not explicitely written in~\cite{CGS}.
However, as the proof of the regularity \Regsoluz\ of the solution 
performed there relies on a~priori estimates and compactness arguments
and the dependence on $\tau$ in the whole calculation of~\cite{CGS} is always made explicit,
\eqref{stab}~holds as well,
and we stress that the corresponding constant $C_0$ does not depend on~$\tau$.

Once well-posedness for problem \State\ is established,
we can address the corresponding control problem.
As in the Introduction, given two functions
\Beq
  \zQ \in \LQ2 
  \aand
  \zS \in \LS2 
  \label{hpzz}
\Eeq
and \pier{three} nonnegative constants $\bQ,\, \bS,\, \bz$,
we~set
\Beq
  \calJ(y,y_\Gamma,\uG)
  := \frac\bQ 2 \, \norma{y-\zQ}_{L^2(Q)}^2
  + \frac\bS 2 \, \norma{y_\Gamma-\zS}_{L^2(\Sigma)}^2
  + \frac\bz 2 \, \norma\uG_{L^2(\Sigma)}^2
  \label{defcost}
\Eeq
for, say, $y\in\L2H$, $\yG\in\L2\HG$ and $\uG\in\LS2$,
and consider the problem of \minimiz ing the cost functional \eqref{defcost}
subject to the constraint $\uG\in\Uad$,
where the control box $\Uad$ is given~by
\Bsist
  & \Uad := 
  & \bigl\{ \uG\in\H1\HG\cap\LS\infty:
  \non
  \\
  && \ \uGmin\leq\uG\leq\uGmax\ \aeS,\ \norma{\dt\uG}_{\LS2}\leq\Mz
  \bigr\} 
  \label{defUad}
\Esist
and to the state system \State.
We simply assume~that
\Beq
  \Mz > 0 , \quad
  \uGmin ,\, \uGmax \in \LS\infty 
  \aand
  \hbox{$\Uad$ is nonempty}.
  \label{hpUad}
\Eeq
Besides, we consider the analogous control problem
of \minimiz ing the cost functional \eqref{defcost}
subject to the constraint $\uG\in\Uad$
and to the state system \Pbltau.
From \cite[Thm.~\pier{2.3}]{CGSvisc}, we~have 
\juerg{the following result.}

\Bthm
\label{Optimumvisc}
Assume \HPstruttura\ and \HPdati,
and let $\calJ$ and $\Uad$ be defined by \eqref{defcost} and \eqref{defUad}
under the assumptions \eqref{hpzz} and~\eqref{hpUad}.
Then, for every $\tau>0$, there exists $\utopt\in\Uad$ such~that
\Beq
  \calJ(\ytopt,\ytoptG,\utopt)
  \leq \calJ(\yt,\ytG,\uG)
  \quad \hbox{for every $\uG\in\Uad$}\,,
  \label{optimumvisc}
\Eeq
where $\ytopt$, $\ytoptG$, $\yt$ and $\ytG$
are the components of the solutions $(\ytopt,\ytoptG,\wtopt)$ and $(\yt,\ytG,\wt)$
to~the state system \Pbltau\
corresponding to the controls $\utopt$ and~$\uG$, respectively.
\Ethm

In~\cite{CGSvisc} \pier{first-order} necessary conditions are obtained 
in terms of the solution to a proper adjoint system.
More precisely, just the case $\tau=1$ is considered there.
However, by going through the paper with some care,
one easily reconstructs the version of the adjoint system corresponding 
to an arbitrary~$\tau>0$.
Even though the adjoint problem considered in~\cite{CGSvisc} 
involves a triplet $(\pt,\qt,\qtG)$ as an adjoint state, 
only the third component $\qtG$ enters the necessary condition for optimality.
On the other hand, $\qt$~and $\qtG$ are strictly related to each other.
Hence, we mention the result that deals with the pair~$(\qt,\qtG)$.
To this end, we recall a tool, the \generaliz ed Neumann problem solver~$\calN$, 
that is often used in connection with the Cahn--Hilliard equations.
With the notation for the mean value introduced in~\eqref{media}, we define 
\Beq
  \dom\calN := \graffe{\vstar\in\Vp: \ \vstar^\Omega = 0}
  \aand
  \calN : \dom\calN \to \graffe{v \in V : \ \vO = 0}
  \label{predefN}
\Eeq
by setting, for $\vstar\in\dom\calN$,
\Beq
  \calN\vstar \in V, \quad
  (\calN\vstar)^\Omega = 0 ,
  \aand
  \iO \nabla\calN\vstar \cdot \nabla z = \< \vstar , z >
  \quad \hbox{for every $z\in V$} .
  \label{defN}
\Eeq
Thus, $\calN\vstar$ is the solution $v$ to the \generaliz ed Neumann problem for $-\Delta$
with datum~$\vstar$ that satisfies~$\vO=0$.
Indeed, if $\vstar\in H$, the above variational equation means \juerg{that}
$-\Delta\calN\vstar = \vstar$ and $\dn\calN\vstar = 0$.
As $\Omega$ is bounded, smooth, and connected,
it turns out that \eqref{defN} yields a well-defined isomorphism.
Moreover, we have
\Beq
  \< \ustar , \calN \vstar >
  = \< \vstar , \calN \ustar >
  = \iO (\nabla\calN\ustar) \cdot (\nabla\calN\vstar)
  \quad \hbox{for $\ustar,\vstar\in\dom\calN$},
  \label{simmN}
\Eeq
whence also
\Beq
  2 \< \dt\vstar(t) , \calN\vstar(t) >
  = \frac d{dt} \iO |\nabla\calN\vstar(t)|^2
  = \frac d{dt} \, \normaVp{\vstar(t)}^2
  \quad \aat
  \label{dtcalN}
\Eeq
for every $\vstar\in\H1\Vp$ satisfying $(\vstar)^\Omega=0$ \aet,
where we have set
\Beq
  \normaVp\vstar := \normaH{\nabla\calN\vstar}
  \quad \hbox{for $\vstar\in\Vp$} .
  \label{defnormaVp}
\Eeq
One easily sees that $\normaVp\cpto$ is a norm in $\Vp$ 
\juerg{which is}
equivalent to the usual dual norm.

Furthermore, we introduce the spaces $\HO$ and $\VO$ \juerg{by setting}
\Beq
  \HO := \graffe{(v,\vG) \in \calH :\ \vO=0}
  \aand
  \VO := \HO \cap \calV\,,
  \label{defHOVO}
\Eeq
and endow them with their natural topologies as subspaces of $\calH$ and~$\calV$, respectively.
As in \cite[Thms.~\pier{2.5} and~5.4]{CGSvisc}, we~have \juerg{the following result.}

\Bthm
\label{Existenceadj}
Assume 
\Beq
  \lam \in \LQ\infty , \quad
  \lamG \in \LS\infty , \quad 
  \phQ \in \LQ2
  \aand
  \phS \in \LS2 .
  \label{genlam}
\Eeq
Then, for every $\tau>0$, there exists a unique pair $(\qt,\qtG)$
satisfying the regularity conditions
\Beq
  \qt \in \H1H \cap \L2{\Hx2}
  \aand
  \qtG \in \H1\HG \cap \L2{\HxG2} 
  \label{regqqG}
\Eeq
and solving the following adjoint problem:
\Bsist
  && (\qt,\qtG)(t) \in \VO
  \quad \hbox{for every $t\in[0,T]$} \vphantom{\frac 1{|\Omega|}}
  \label{primaN}
  \\
  && - \iO \dt \bigl( \calN(\qt(t)) + \tau \qt(t) \bigr) \, v
  + \iO \nabla\qt(t) \cdot \nabla v
  + \iO \lam(t) \, \qt(t) \, v
  \qquad
  \non
  \\
  && \qquad {}
  - \iG \dt\qtG(t) \, \vG
  + \iG \nablaG\qtG(t) \cdot \nablaG\vG
  + \iG \lamG(t) \, \qtG(t) \, \vG
  \non
  \\
  && = \iO \phQ(t) v
  + \iG \phS(t) \vG
  \quad \hbox{\aat\ and every $(v,\vG)\in\VO$}
  \qquad
  \label{secondaN}
  \\
  && \iO \bigl( \calN\qt + \tau \qt \bigr)(T) \, v 
  + \iG \qG(T) \, \vG
  = 0
  \quad \hbox{for every $(v,\vG)\in\VO$} \,.
  \label{cauchyN}
\Esist
\Accorpa\AggiuntoN primaN cauchyN
\Ethm
More precisely, in \cite{CGSvisc} 
the above theorem is proved with the particular choice
\Beq
  \lam = f''(\ytopt) , \quad
  \lamG = \fG''(\ytoptG) , \quad
  \phQ = \bQ (\ytopt-\zQ) 
  \aand
  \phS = \bS (\ytoptG-\zS) 
  \label{scelteCGS}
\Eeq
where $\ytopt$ and $\ytoptG$ are the components of the state 
associated to an optimal control~$\utopt$.
However, the same proof is valid under assumption~\eqref{genlam}.

Finally, \cite[Thm.~\pier{2.6}]{CGSvisc} gives a necessary condition 
for $\utopt$ to be an optimal control
in terms of the solution to the above adjoint system
corresponding to~\eqref{scelteCGS}.
\juerg{This} condition reads
\Beq
  \intS (\qtG + \bz \utopt) (\vG - \utopt) \geq 0
  \quad \hbox{for every $\vG\in\Uad$}.
  \label{cnoptAV}
\Eeq

In this paper, we first show the existence of an optimal control~$\uopt$.
Namely, we prove the following result.

\Bthm
\label{Optimum}
Assume \HPstruttura\ and \HPdati,
and let $\calJ$ and $\Uad$ be defined by \eqref{defcost} and \eqref{defUad}
under the assumptions \eqref{hpzz} and~\eqref{hpUad}.
Then there exists \juerg{some} $\uopt\in\Uad$ such~that
\Beq
  \calJ(\yopt,\yoptG,\uopt)
  \leq \calJ(y,\yG,\uG)
  \quad \hbox{for every $\uG\in\Uad$}\,,
  \label{optimum}
\Eeq
where $\yopt$, $\yoptG$, $y$ and $\yG$
are the components of the solutions $(\yopt,\yoptG,\wopt)$ and $(y,\yG,w)$
to~the state system \State\
corresponding to the controls $\uopt$ and~$\uG$, respectively.
\Ethm

Next, for every optimal control~$\uopt$,
we derive a necessary optimality condition like \eqref{cnoptAV}
in terms of the solution of a \generaliz ed adjoint system.
In order to make the last sentence precise,
we introduce the spaces
\Bsist
  && \calW := \L2\VO \cap \bigl( \H1\Vp \times \H1\VGp \bigr)
  \label{defZ}
  \\
  && \calWz := \graffe{ (v,\vG) \in \calW :\ (v,\vG)(0) = (0,0) }
  \qquad
  \label{defZz}
\Esist
and endow them with their natural topologies.
Moreover, we denote by $\[\cpto,\cpto]$
the duality product between $\calWp$ and~$\calWz$.
We will prove the following representation result for the elements 
of the dual space~$\calWp$.

\Bprop
\label{ReprWp}
A functional $F:\calWz\to\erre$ belongs to $\calWp$ if and only if
there exist $\Lam$ and $\LamG$ satisfying
\Bsist
  && \Lam \in \bigl( \H1\Vp \cap \L2V \bigr)^*
  \aand
  \LamG \in \bigl( \H1\VGp \cap \L2\VG \bigr)^*
  \qquad\quad
  \label{goodLam}
  \\
  && \[ F , (v,\vG) ]
  = \< \Lam , v >_Q
  + \< \LamG , \vG >_\Sigma
  \quad \hbox{for every $(v,\vG)\in\calWz$}\,, 
  \label{reprWp}
\Esist
where the duality products $\<\cpto,\cpto>_Q$ and $\<\cpto,\cpto>_\Sigma$
are related to the spaces $X^*$ and $X$ with
$X=\H1\Vp \cap \L2V$ and $X=\H1\VGp \cap \L2\VG$, respectively.
\Eprop

However, \juerg{this} representation is not unique,
since different pairs $(\Lam,\LamG)$ satisfying \eqref{goodLam}
could \juerg{generate} the same functional $F$ through formula~\eqref{reprWp}.

At this point, we are ready to present our result 
on the necessary optimality conditions 
for the control problem related to the pure Cahn--Hilliard equations, i.e.,
the analogue of \eqref{cnoptAV}
in terms of a solution to a \generaliz ed adjoint system.

\Bthm
\label{CNopt}
Assume \HPstruttura\ and \HPdati,
and let $\calJ$ and $\Uad$ be defined by \eqref{defcost} and \eqref{defUad}
under the assumptions \eqref{hpzz} and~\eqref{hpUad}.
Moreover, let $\uopt$ be any optimal control 
as in the statement of Theorem~\ref{Optimum}.
Then, there exist $\Lam$ and $\LamG$ satisfying \eqref{goodLam},
and a pair $(q,\qG)$ satisfying
\Bsist
  && q \in \L\infty\Vp \cap \L2V
  \label{regq}
  \\
  && \qG \in \L\infty\HG \cap \L2\VG
  \label{regqG}
  \\
  && (q,\qG)(t) \in \VO
  \quad \hbox{for every $t\in[0,T]$}\,,
  \label{primazero}
\Esist
as well as 
\Bsist
  && \ioT \< \dt v(t) , \calN q(t) > \, dt
  + \ioT \< \dt\vG(t) , \qG(t) >_\Gamma \, dt
  \non
  \\
  && \quad {}
  + \intQ \nabla q \cdot \nabla v
  + \intS \nablaG\qG \cdot \nablaG\vG
  + \< \Lam , v >_Q
  + \< \LamG , \vG >_\Sigma
  \qquad
  \non
  \\
  && = \intQ \phQ \, v
  + \intS \phS \, \vG
  \qquad \hbox{for every $(v,\vG)\in\calWz$}\,, 
  \label{secondazero}
\Esist
such that
\Beq
  \intS (\qG + \bz \uopt) (\vG - \uopt) \geq 0
  \quad \hbox{for every $\vG\in\Uad$}.
  \label{cnopt}
\Eeq
In particular, if $\bz>0$, \juerg{then} the optimal control $\uopt$
is the $\LS2$-projection of $-\qG/\bz$ \juerg{onto} $\Uad$.
\Ethm

One \recogniz es in \eqref{secondazero}
a problem that is analogous to \accorpa{secondaN}{cauchyN}.
Indeed, if $\Lam$, $\LamG$ and the solution $(q,\qG)$ were regular functions,
\juerg{then its strong form} should contain 
both a \generaliz ed \pier{backward} parabolic equation like \eqref{secondaN} 
and a final condition for $(\calN q,\qG)$ of type \eqref{cauchyN},
since the definition of $\calWz$ allows its elements to be free at $t=T$.
However, the terms $\lamt\qt$ and $\lamtG\qtG$ are just replaced
by the functionals $\Lam$ and~$\LamG$ and 
\juerg{cannot be} identified as products, unfortunately.


\section{Proofs}
\label{PROOFS}
\setcounter{equation}{0}

In the whole section, we assume that all \juerg{of} the conditions
\HPstruttura\ and \accorpa{hpyz}{hpyzbis}
on the structure and the initial datum of the state system,
as well as assumptions \eqref{hpzz} and~\eqref{hpUad} 
that regard the cost functional \eqref{defcost} and the control box \eqref{defUad},
are satisfied.
We start with an expected result.

\Bprop
\label{Convergence}
Assume $\utG\in\H1\HG$ and let $(\yt,\ytG,\wt)$
be the solution to \pier{the} problem \Pbltau\ associated to~$\utG$.
If $\utG$ converges to $\uG$ \pier{weakly} in $\H1\HG$ as $\tau\seto0$, then
\Bsist
  & \yt \to y
  & \quad \hbox{weakly star in $\H1\Vp\cap\L\infty V\cap\L2\Hdue$}
  \non
  \\
  && \qquad \hbox{and strongly in \pier{$\C0H \cap \L2V$}}
  \label{convyt}
  \\
  & \ytG \to \yG
  & \quad \hbox{weakly star in $\H1\HG\cap\L\infty\VG\cap\L2\HdueG$}
  \non
  \\
  && \qquad \hbox{and strongly in \pier{$\C0\HG \cap \L2\VG$}}
  \label{convytG}
  \\
  & \wt \to w
  & \quad \hbox{weakly star in $\pier{\L2V}$}\,,
  \label{convwt}
\Esist
where $(y,\yG,w)$ is the solution to problem \State\
associated \juerg{with} $\uG$.
\Eprop

\Bdim
The family $\graffe{\utG}$ is bounded in $\H1\HG$.
Thus, the solution $(\yt,\ytG,\wt)$ satisfies~\eqref{stab}
for some constant~$C_0$, so that the weak or weak star convergence
specified in \accorpa{convyt}{convwt} holds for a subsequence.
In particular, the Cauchy condition \eqref{cauchy} for $y$ is satisfied.
Moreover, we also have
$\tau\,\dt\yt\to0$ strongly in $\L2H$ as well as
$f'(\yt)\to\xi$ and $\fG'(\ytG)\to\xiG$
weakly in $\L2H$ and in $\L2\HG$, respectively,
for some $\xi$ and $\xiG$.
Furthermore, $\yt$~and $\ytG$ converge to their limits strongly
in $\L2H$ and $\L2\HG$, respectively, thanks to the Aubin-Lions lemma
(see, e.g., \cite[Thm.~5.1, p.~58]{Lions}, which also implies a much better strong convergence \cite[Sect.~8, Cor.~4]{Simon}).
Now, as said in Section~\ref{STATEMENT}, 
we can split $f'$ as $f'=\beta+\pi$, 
where $\beta$ is monotone \juerg{and} $\pi$ is \Lip\ continuous.
It follows that $\pi(\yt)$ converges to $\pi(y)$ strongly in $\L2H$,
whence \juerg{we obtain that also} $\beta(\yt)$ converges to $\xi-\pi(y)$ weakly in~$\L2H$.
Then, we infer that $\xi-\pi(y)=\beta(y)$ \aeQ, i.e., $\xi=f'(y)$ \aeQ,
with the help of standard monotonicity arguments
(see, e.g., \cite[Lemma~1.3, p.~42]{Barbu}).
Similarly, we have $\xiG=\fG'(\yG)$.
Therefore, by starting from \eqref{intprima} and \eqref{intsecondavisc} 
written with $\utG$ in place of~$\uG$,
we can pass to the limit and obtain \IntPbl\
associated to the limit control~$\uG$.
As the solution to the limit problem is unique,
the whole family $(\yt,\ytG,\wt)$ converges to $(y,\yG,w)$ in the sense of the statement
and the proof is complete.
\Edim

\Bcor
\label{CorConvergence}
Estimate \eqref{stab}, \juerg{written formally} with $\tau=0$,
holds for the solution to the pure Cahn--Hilliard system \State.
\Ecor

\Bdim
By applying the above proposition with $\utG=\uG$
and using \eqref{stab} for the solution to the viscous problem,
we immediately conclude \juerg{the claim}.
\Edim

\step
Proof of Theorem~\ref{Optimum}

We use the direct method and start from a \minimiz ing sequence $\graffe{\unG}$.
Then, $\unG$~remains bounded in $\H1H$, whence we have $\unG\to\uopt$ weakly in $\H1\HG$ for a subsequence.
By Corollary~\ref{CorConvergence},
the sequence of the corresponding states $(\yn,\ynG,\wn)$ satisfies the analogue of~\eqref{stab}.
Hence, by arguing as in the proof of Proposition~\ref{Convergence},
we infer that the solutions $(\yn,\ynG,\wn)$ converge in the proper topology
to the solution $(\yopt,\yoptG,\wopt)$ associated to~$\uopt$.
In particular, there holds the strong convergence 
specified by the analogues of \eqref{convyt} and~\eqref{convytG}.
Thus, by also owing to \juerg{the semicontinuity of $\calJ$ and the optimality of}
$\unG$, we have
\Beq
  \calJ(\yopt,\yoptG,\uopt) \leq \liminf_{n\to\infty} \calJ(\yn,\ynG,\unG)
  \leq \calJ(y,\yG,\uG)
  \non
\Eeq
for every $\uG\in\Uad$, where $y$ and $\yG$ are the components 
of the solution to the Cahn--Hilliard system associated \juerg{with} $\uG$.
This means that $\uopt$ is an optimal control.\qed

\step
Proof of Proposition~\ref{ReprWp}

Assume \pier{that $\Lam$ and $\LamG$} satisfy~\eqref{goodLam}.
Then, formula \eqref{reprWp} actually defines a functional $F$ on~$\calWz$.
Clearly, $F$~is linear.
Moreover, we have, for every $(v,\vG)\in\calWz$,
\Bsist
  && |\<\Lam,v>_Q + \<\LamG,\vG>_\Sigma|
  \non
  \\
  && \leq \norma\Lam_{(\H1\Vp\cap\L2V)^*} \, \norma v_{\H1\Vp\cap\L2V}
  \non
  \\
  && \quad {}
  + \norma\LamG_{(\H1\VGp\cap\L2\VG)^*} \, \norma v_{\H1\VGp\cap\L2\VG}
  \non
  \\
  && \leq \bigl( \norma\Lam_{(\H1\Vp\cap\L2V)^*} + \norma\LamG_{(\H1\VGp\cap\L2\VG)^*} \bigr) \norma{(v,\vG)}_{\calW}\,,
  \non
\Esist
so that $F$~is continuous.
Conversely, assume that $F\in\calWp$.
As $\calWz$ is~a (closed) subspace of $\tilde\calW:=\bigl(\H1\Vp\cap\L2V)\times(\H1\VGp\cap\L2\VG)$,
we can extend $F$ to a linear continuous functional $\tilde F$ on~$\tilde\calW$.
Then, there exist $\Lam$ and~$\LamG$ 
(take $\Lam(v):=\tilde F(v,0)$ and $\LamG(\vG):=\tilde F(0,\vG)$)
satisfying \eqref{goodLam} such that
\Beq
  \< \tilde F , (v,\vG) >
  = \< \Lam , v >_Q
  + \< \LamG , \vG >_\Sigma
  \quad \hbox{for every $(v,\vG)\in\tilde\calW$}\,,
  \non
\Eeq
where the duality product on the \lhs\ refers to the spaces $(\tilde\calW)^*$ and~$\tilde\calW$.
Since $\[F,(v,\vG)]=\<\tilde F,(v,\vG)>$ for every $(v,\vG)\in\calWz$,
\eqref{reprWp} immediately follows.\qed

\bigskip

The rest of this section is devoted to the proof of Theorem~\ref{CNopt}
on the necessary optimality conditions.
Therefore, \pier{besides} the general assumptions, we also suppose~that
\Beq
  \hbox{\sl $\uopt$ is any optimal control as in Theorem~\ref{Optimum}},
  \label{uoptfixed}
\Eeq
that is, an arbitrary optimal control $\uopt$ is fixed once and for all.
In order to arrive at the desired necessary optimality condition for~$\uopt$,
we follow~\cite{Barbutrick} and introduce the \pier{modified} functional $\tilde\calJ$ defined \pier{by}
\Beq
  \tilde\calJ(y,\yG,\uG)
  := \calJ(y,\yG,\uG) + \frac 12 \, \norma{\uG - \uopt}_{\LS2}^2 \,.
  \label{modifiedcost}
\Eeq
Then the analogue of Theorem~\ref{Optimumvisc} holds,
and we have:

\Bthm
\label{Optimummod}
For every $\tau>0$, there exists \juerg{some} $\utmod\in\Uad$ such~that
\Beq
  \tilde\calJ(\ytmod,\ytmodG,\utmod)
  \leq \tilde\calJ(\yt,\ytG,\uG)
  \quad \hbox{for every $\uG\in\Uad$}\,,
  \label{optimummod}
\Eeq
where $\ytmod$, $\ytmodG$, $\yt$ and $\ytG$
are the components of the solutions $(\ytmod,\ytmodG,\wtmod)$ and $(\yt,\ytG,\wt)$
to~the state system \Pbltau\
corresponding to the controls $\utmod$ and~$\uG$, respectively.
\Ethm

For the reader's convenience, we fix the notation just used
and introduce a new one
(\juerg{which was} already used with a different meaning earlier in this paper):
\Bsist
  && \hbox{\sl $\utmod$ is an optimal control as in Theorem~\ref{Optimummod}}
  \label{utmodfixed}
  \\
  && \hbox{\sl $(\ytmod,\ytmodG,\wtmod)$ is the solution to \Pbltau\ corresponding to $\utmod$}
  \qquad
  \label{taustate}
  \\
  && \hbox{\sl $(\ytopt,\ytoptG,\wtopt)$ is the solution to \Pbltau\ corresponding to $\uopt$}.
  \qquad
  \label{newstate}
\Esist
The next step consists in writing the proper adjoint system 
and the corresponding necessary optimality condition, 
\juerg{which} can be done by repeating the argument of~\cite{CGSvisc}.
However, instead of just stating the corresponding result,
we spend some words that can help the reader.
The optimality variational inequality is derived as a condition
on the \Frechet\ derivative of the map
(defined in a proper functional framework)
$\uG\mapsto\tilde\calJ(y,\yG,\uG)$,
where the pair $(y,\yG)$ is subjected to the state system.
Thus, \juerg{this} derivative depends
on the \Frechet\ derivative of the functional
$(y,\yG,\uG)\mapsto\tilde\calJ(y,\uG,\uG)$,
which is given~by
\Beq
  [D\tilde\calJ(y,\pier{\yG},\uG)](k,k_\Gamma,\hG)]
  = \bQ \intQ (y-\zQ) k
  + \bS \intS (\yG-\zS) k_\Gamma
  + \intS \bigl( \bz\uG + (\uG-\uopt) \bigr) \hG \,.
  \non
\Eeq
Hence, the argument for $\tilde\calJ$ differs from the one for $\calJ$
only in relation to the last integral.
In other words, we just have to replace
$\bz\uG$ by $\bz\uG+(\uG-\uopt)$ in the whole argument of~\cite{CGSvisc}.
In particular, the adjoint system remains unchanged.
Here is the conclusion.

\Bprop
\label{CNoptmod}
With the notations \accorpa{utmodfixed}{taustate},
we have
\Beq
  \intS \bigl( \qtG + \bz\utmod + (\utmod - \uopt) \bigr) (\vG - \utmod) \geq 0
  \quad \hbox{for every $\vG\in\Uad$}\,,
  \label{cnoptmod}
\Eeq
where $\qtG$ is the component of the solution $(\qt,\qtG)$ to
\AggiuntoN\ corresponding \pier{to} $\uG=\utmod$ with the 
choices $\lam=\lamt$, $\lamG=\lamtG$, $\phQ=\phtQ$ and $\phS=\phtS$ specified~by
\Beq
  \lamt=f''(\ytmod), \enskip
  \lamtG=\fG''(\ytmodG), \enskip 
  \phtQ=\bQ(\ytmod-\zQ) \enskip and \enskip
  \phtS=\bS(\ytmodG-\zS) .
  \label{sceltelam}
\Eeq
\Eprop

Thus, our project for the proof of Theorem~\ref{CNopt} is the following:
we take the limit in \eqref{cnoptmod} 
and in the adjoint system mentioned in the previous statement as $\tau$ tends to zero.
This will lead to the desired necessary optimality condition \eqref{cnopt}
provided that we prove that the optimal controls $\utmod$ converge to~$\uopt$.
The \juerg{details of this project are} the following.

\noindent
$i)$~There hold
\Bsist
  & \utmod \to \uopt
  & \quad \hbox{weakly \pier{star in $\H1\HG\cap L^\infty (\Sigma)$} and strongly in $\LS2$} \qquad
  \label{convu}
  \\
  \separa
  & \ytmod \to \yopt
  & \quad \hbox{weakly star in $\H1\Vp\cap\L\infty V\cap\L2\Hdue$}
  \non
  \\
  && \qquad \hbox{and strongly in \pier{$\C0H \cap \L2V$}}
  \label{convy}
  \\
  \separa
  & \ytmodG \to \yoptG
  & \quad \hbox{weakly star in $\H1\HG\cap\L\infty\VG\cap\L2\HdueG$}
  \qquad
  \non
  \\
  && \qquad \hbox{and strongly in \pier{$\C0\HG \cap \L2\VG$}}
  \label{convyG}
  \\
  & \wtmod \to \wopt
  & \quad \hbox{weakly star in $\pier{\L2V}$} 
  \label{convw}
  \\
  & \qt \to q
  & \quad \hbox{weakly star in $\L\infty\Vp\cap\L2V$}
  \label{convq}
  \\
  & \qtG \to \qG
  & \quad \hbox{weakly star in $\L\infty\HG\cap\L2\VG$}\,,
  \qquad\qquad
  \label{convqG}
\Esist
as well as
\Beq
  \tilde\calJ(\ytmod,\ytmodG,\utmod)
    \to \calJ(\yopt,\yoptG,\uopt)\,,
  \label{convJ}
\Eeq
at least for a subsequence,
and $(\yopt,\yoptG,\wopt)$ solves problem \State\ with $\uG=\uopt$.

\noindent
$ii)$~The functionals associated \juerg{with} the pair $(\lamt\qt,\lamtG,\qtG)$ 
through Proposition~\ref{ReprWp}
converge to some functional weakly in~$\calWp$, at least for a subsequence,
and we then represent the limit by some pair~$(\Lam,\LamG)$,
so that we have
\Beq
  \<\lamt\qt,v>_Q + \<\lamtG\qtG,\vG>_\Sigma
  \to \<\Lam,v>_Q + \<\LamG,\vG>_\Sigma
  \quad \hbox{for every $(v,\vG)\in\calWz$}.
  \label{convlambda}
\Eeq

\noindent
$iii)$~With such a choice of $(\Lam,\LamG)$, 
the pair $(q,\qG)$ solves \accorpa{primazero}{secondazero}.

\noindent
$iv)$~Condition \eqref{cnopt} holds.

The main tool is proving a priori estimates.
To this concern, we use the following rule to denote constants
in order to avoid a boring notation.
The small-case symbol $c$ stands for different constants 
\juerg{that neither} depend on $\tau$ nor on the functions whose norm we want to estimate.
Hence, the meaning of $c$ might
change from line to line and even in the same chain of equalities or inequalities.
Similarly, a~symbol like $c_\delta$ denotes different constants
that depend on the parameter~$\delta$, in addition.

\step
First a priori estimate

As $\utopt\in\Uad$ and Theorem~\ref{daCGS} holds, we have
\Bsist
  \hskip-1cm && \norma\utmod_{\H1\HG}
  + \norma\ytmod_{\H1\Vp \cap \L\infty V \cap \L2\Hdue}
  \non
  \\[1mm]
  \hskip-1cm && \quad {}
  + \norma\ytmodG_{\H1\HG \cap \L\infty\VG \cap \L2\HdueG}
  + \norma\wtmod_{\pier{\L2V}}
  \non
  \\[1mm]
  \hskip-1cm && \quad {}
  + \norma{f'(\ytmod)}_{\L2H}
  + \norma{\fG'(\ytmodG)}_{\L2\HG}
  + \tau^{1/2} \norma{\dt\ytmod}_{\L2H}
  \leq c \,.
  \label{primastima}
\Esist

\step
Second a priori estimate

For the reader's convenience, we 
\juerg{explicitly} write the adjoint system mentioned in Proposition~\ref{CNoptmod},
as well as the regularity of its solution,
\Bsist
  && \qt \in \H1H \cap \L2{\Hx2} , \quad
  \qtG \in \H1\HG \cap \L2{\HxG2} 
  \qquad
  \label{regqGbis}
  \\
  && (\qt,\qtG)(s) \in \VO
  \quad \hbox{for every $s\in[0,T]$} \vphantom{\frac 1{|\Omega|}}
  \label{primaNbis}
  \\
  \separa
  && - \iO \dt \bigl( \calN(\qt(s)) + \tau \qt(s) \bigr) \, v
  + \iO \nabla\qt(s) \cdot \nabla v
  + \iO \lamt(s) \, \qt(s) \, v
  \qquad
  \non
  \\
  && \qquad {}
  - \iG \dt\qtG(s) \, \vG
  + \iG \nablaG\qtG(s) \cdot \nablaG\vG
  + \iG \lamtG(s) \, \qtG(s) \, \vG
  \non
  \\
  && = \iO \phtQ(s) v
  + \iG \phtS(s) \vG
  \quad \hbox{for a.e.\ $s\in(0,T)$ and every $(v,\vG)\in\VO$}
  \qquad
  \label{secondaNbis}
  \\
  \separa
  && \quad \hbox{where
     $\lamt=f''(\ytmod)$, \
     $\lamtG=\fG''(\ytmodG)$, \ 
     $\phtQ=\bQ(\ytmod-\zQ)$ \ and \
     $\phtS=\bS(\ytmodG-\zS)$
     }\vphantom\int
  \non
  \\
  && \iO \bigl( \calN\qt + \tau \qt \bigr)(T) \, v 
  + \iG \qG(T) \, \vG
  = 0
  \quad \hbox{for every $(v,\vG)\in\VO$} .
  \label{cauchyNbis}
\Esist
Now, we choose $v=\qt(s)$ and $\vG=\qtG(s)$,
and integrate over $(t,T)$ with respect to~$s$.
\pier{Recalling} \eqref{dtcalN} and now reading
$Q_t:=\Omega\times(t,T)$ and $\Sigma_t:=\Gamma\times(t,T)$,
we~have
\Bsist
  && \frac 12 \, \normaVp{\qt(t)}^2
  + \frac \tau 2 \iO |\qt(t)|^2
  + \intQt |\nabla\qt|^2
  + \intQt \lamt |\qt|^2
  \non
  \\
  && \quad {}
  + \frac 12 \iG |\qtG(t)|^2
  + \intSt |\nablaG\qtG|^2
  + \intSt \lamtG |\qtG|^2
  \non
  \\
  && = \intQt \phtQ\, \qt
  + \intSt \phtS \, \qtG
  \leq \intQ |\phtQ|^2 + \intQt |\qt|^2
  + \intS |\phtS|^2 + \intSt |\qtG|^2 
  \non
  \\
  && \leq \intQt |\qt|^2
  + \intSt |\qtG|^2 
  + c 
  \label{persecondastima}
\Esist
\juerg{where the last inequality follows from} \eqref{primastima}.
By accounting for~\eqref{hpfseconda}, we also have
\Beq
  \intQt \lamt |\qt|^2
  \geq - c \intQt |\qt|^2
  \aand
  \intQt \lamtG |\qtG|^2
  \geq - c \intSt |\qtG|^2 .
  \non
\Eeq
We treat the volume integral (and the same on the \rhs\ of~\eqref{persecondastima})
invoking the compact embedding $V\subset H$.
We have
\Beq
  \iO |v|^2
  \leq \delta \iO |\nabla v|^2 + c_\delta \normaVp v^2
  \quad \hbox{for every $v\in V$ and $\delta>0$}.
  \non
\Eeq
Hence, we deduce \juerg{that}
\Beq
  \intQt |\qt|^2
  \leq \delta \intQt |\nabla\qt|^2 + c_\delta \int_t^T \normaVp{\qt(s)}^2 \, ds \,.
  \non
\Eeq
Therefore, by combining, choosing $\delta$ small enough and applying the backward Gronwall lemma, 
we conclude~that
\Beq
  \norma\qt_{\L\infty\Vp\cap\L2V}
  + \norma\qtG_{\L\infty\HG\cap\L2\VG}
  + \tau^{1/2} \norma\qt_{\L\infty H}
  \leq c \,.
  \label{secondastima}
\Eeq

\step
Third a priori estimate

\pier{Take an arbitrary pair} $(v,\vG)\in\H1\HO\cap\L2\VO$,
\pier{and} test \eqref{secondaNbis} by $v(s)$ and $\vG(s)$.
Then, we sum over $s\in(0,T)$ and integrate by parts
with the help of~\eqref{cauchyNbis}, 
so that no integral related to the time $T$ appears.
In particular, if $(v,\vG)\in\calWz$, even the terms 
evaluated at $t=0$ vanish
and we obtain \juerg{that}
\Bsist
  && \intQ (\calN\qt+\tau\qt) \dt v
  + \intQ \nabla\qt \cdot \nabla v
  + \intQ \lamt \qt v
  + \intS \qtG \dt\vG
  + \intS \nabla\qtG \cdot \nabla\vG
  + \intS \lamt \qtG \vG
  \non
  \\
  && = \intQ \phtQ \, v
  + \intS \phtS \, \vG \,.
  \label{adjperparti}
\Esist
Therefore, we have, for every $(v,\vG)\in\calWz$,
\Bsist
  && \left| \intQ \lamt \qt v + \intS \lamt \qtG \vG \right|
  \non
  \\
  && \leq \norma{\calN\qt+\tau\qt}_{\L2V} \, \norma{\dt v}_{\L2\Vp} 
  + \norma\qt_{\L2V} \, \norma v_{\L2V}
  \non
  \\
  && \quad {}
  + \norma\qtG_{\L2\VG} \, \norma{\dt \pier{\vG}}_{\pier{\L2\VGp}} 
  + \norma{\pier{\qtG}}_{\pier{\L2\VG}} \, \norma\vG_{\L2\VG}
  \non
  \\
  && \quad {}
  + \norma\phtQ_{\L2H} \, \norma v_{\L2H}
  + \norma\phtS_{\L2\HG} \, \norma\vG_{\L2\HG} \,.
  \non
\Esist
Now, by assuming $\tau\leq1$, we have
$\normaV{\calN v+\tau v}\leq c\normaVp v+\tau\normaV v\leq c\normaV v$ 
for every $v\in V$ with zero mean value 
(see~\eqref{defnormaVp}).
Therefore, by accounting for \eqref{primastima} and \eqref{secondastima},
we conclude~that
\Beq
  \left| \intQ \lamt \qt v + \intS \lamt \qtG \vG \right|
  \leq c \, \norma{(v,\vG)}_\calW
  \quad \hbox{for every $(v,\vG)\in\calWz$}.
  \label{terzastima}
\Eeq

\step
Conclusion of the proof of Theorem~\ref{CNopt}

From the above estimates, we infer that
\Bsist
  & \utmod \to \uG
  & \quad \hbox{weakly \pier{star in $\H1\HG\cap L^\infty (\Sigma)$}}
  \label{perconvu}
  \\
  \separa
  & \ytmod \to y
  & \quad \hbox{weakly star in $\H1\Vp\cap\L\infty V\cap\L2\Hdue$}\qquad
  \non
  \\
  && \qquad \hbox{and strongly in \pier{$\C0H \cap \L2V$}}
  \label{perconvy}
  \\
  \separa
  & \ytmodG \to \yG
  & \quad \hbox{weakly star in $\H1\HG\cap\L\infty\VG\cap\L2\HdueG$}
  \qquad
  \non
  \\
  && \qquad \hbox{and strongly in \pier{$\C0\HG \cap \L2\VG$}}
  \label{perconvyG}
  \\
  & \wtmod \to w
  & \quad \hbox{weakly star in $\pier{\L2V}$} 
  \label{perconvw}
  \\
  & \qt \to q
  & \quad \hbox{weakly star in $\L\infty\Vp\cap\L2V$}
  \label{perconvq}
  \\
  & \qtG \to \qG
  & \quad \hbox{weakly star in $\L\infty\HG\cap\L2\VG$}
  \qquad\qquad
  \label{perconvqG}
  \\
  & \tau \qt \to 0
  & \quad \hbox{\pier{strongly} in $\L\infty H$}
  \label{pertauqt}
\Esist
at least for a subsequence,
and $(y,\yG,w)$ is the solution to \pier{the} problem \Pbl\ 
corresponding to $\uG$,
thanks to Proposition~\ref{Convergence}.
Notice that \accorpa{perconvq}{perconvqG} coincide with \accorpa{convq}{convqG}
and that \accorpa{convu}{convw} hold once we prove that $\uG=\uopt$
and that $\utmod$ \pier{converges} \juerg{strongly in~$\LS2$}.

\juerg{To this end, we} recall the notations \accorpa{utmodfixed}{newstate},
and it is understood that all the limits we write are referred to the selected subsequence.
By optimality, we have
\Beq
  \calJ(\yopt,\yoptG,\uopt) \leq \calJ(y,\yG,\uG)
  \aand
  \tilde\calJ(\ytmod,\ytmodG,\utmod) \leq \tilde\calJ(\ytopt,\ytoptG,\uopt) .
  \non
\Eeq
On the other hand, \accorpa{perconvu}{perconvyG} and Proposition~\ref{Convergence}
applied with $\utG=\uopt$ yield
\Beq
  \tilde\calJ(y,\yG,\uG) \leq \liminf \tilde\calJ(\ytmod,\ytmodG,\utmod)
  \aand
  \lim \calJ(\ytopt,\ytoptG,\uopt) = \calJ(\yopt,\yoptG,\uopt) .
  \non
\Eeq
By combining, we deduce that
\Bsist
  && \calJ(\yopt,\yoptG,\uopt) + \frac 12 \, \norma{\uG-\uopt}_{\LS2}^2
  \leq \calJ(y,\yG,\uG) + \frac 12 \, \norma{\uG-\uopt}_{\LS2}^2
  \non
  \\
  && = \tilde\calJ(y,\yG,\uG)
  \leq \liminf \tilde\calJ(\ytmod,\ytmodG,\utmod)
  \leq \limsup \tilde\calJ(\ytmod,\ytmodG,\utmod)
  \vphantom{\frac 12}
  \non
  \\
  && \leq \limsup \tilde\calJ(\ytopt,\ytoptG,\uopt)
  = \limsup \calJ(\ytopt,\ytoptG,\uopt)
  = \calJ(\yopt,\yoptG,\uopt) .
  \non
\Esist
By comparing the first and last terms of this chain,
we infer that the $\LS2$-norm of $\uG-\uopt$ vanishes,
whence $\uG=\uopt$, as desired.
In order to prove the strong convergence mentioned in~\eqref{convu},
we observe that the above argument also shows~that
\Beq
  \liminf \tilde\calJ(\ytmod,\ytmodG,\utmod)
  = \limsup \tilde\calJ(\ytmod,\ytmodG,\utmod)
  = \calJ(\yopt,\yoptG,\uopt).
  \non
\Eeq
Notice that this coincides with \eqref{convJ}.
From the strong convergence given by \eqref{convy} and \eqref{convyG},
and by comparison, we deduce~that
\Beq
  \lim \left(
    \frac \bz 2 \intS |\utmod|^2 + \frac 12 \intS |\utmod-\uopt|^2
  \right)
  = \frac \bz 2 \intS |\pier{\uopt}|^2\,,
  \non
\Eeq
whence also
\Bsist
  && \limsup \frac \bz 2 \intS |\utmod|^2
  \leq \limsup \left(
    \frac \bz 2 \intS |\utmod|^2 + \frac 12 \intS |\utmod-\uopt|^2
  \right)
  \non
  \\
  && = \frac \bz 2 \intS |\uopt|^2
  \leq \liminf \frac \bz 2 \intS |\utmod|^2 .
  \non
\Esist
Therefore\pier{, we have}
\Beq
  \lim \frac \bz 2 \intS |\utmod|^2
  = \frac \bz 2 \intS |\uopt|^2\,,
  \quad \hbox{whence} \quad
  \lim \frac 12 \intS |\utmod-\uopt|^2 = 0 \,, 
  \non
\Eeq
and \accorpa{convu}{convw} are completely proved.

Now, we deal with the limit $(q,\qG)$ given 
by~\accorpa{perconvq}{perconvqG}, i.e., \accorpa{convq}{convqG}.
Clearly, \eqref{primazero} holds as well.
Furthermore, as $\normaV{\calN\vstar}\leq c\normaVp\vstar$ 
for every $\vstar\in\Vp$ with zero mean value
(see~\eqref{defN}),
\juerg{and since the} convergence \eqref{pertauqt} holds,
we~also have
\Beq
  \calN\qt + \tau \qt \to \calN q
  \quad \hbox{weakly star in $\L\infty H$}.
  \non
\Eeq
On the other hand, \eqref{terzastima} implies that
the functionals $F^\tau\in\calWp$ defined~by
\Beq
  \[ F^\tau , (v,\vG) ]
  := \< \lamt \qt , v >_Q
  + \< \lamtG \qtG , \vG >_\Sigma\,,
  \non
\Eeq
i.e., the functionals associated \juerg{with} $(\lamt\qt,\lamtG,\qtG)$
as in Proposition~\ref{ReprWp}, are bounded in~$\calWp$.
Therefore, for a subsequence, we have $F^\tau\to F$ weakly star in~$\calWp$,
where $F$ is some element of~$\calWp$.
Hence, if we represent $F$ as stated in Proposition~\ref{ReprWp},
we find $\Lam$ and $\LamG$ satisfying \eqref{goodLam} and~\eqref{convlambda}.
At this point, it is \sfw\ to pass to the limit in~\eqref{adjperparti} and in \eqref{cnoptmod}
\juerg{to} obtain both \eqref{secondazero} and~\eqref{cnopt}.
This completes the proof of Theorem~\ref{CNopt}.\qed

\Brem
\label{Fullconv}
The above proof can be repeated without any change
starting from any sequence $\tau_n\seto0$.
By doing that, we obtain:
there exists a subsequence $\graffe{\tau_{n_k}}$ such that
\accorpa{convu}{convJ} hold along the selected subsequence.
As the limits $\uopt$, $\yopt$, $\yoptG$, $\wopt$ and $\calJ(\yopt,\yoptG,\uopt)$
are always the same,
this proves that in fact \accorpa{convu}{convw} as well as \eqref{convJ}
hold for the whole family.
On the contrary, the limits $q$ and $\qG$ might depend on the selected subsequence
since no uniqueness result for the adjoint problem is known.
Nevertheless, the necessary optimality condition \eqref{cnopt}
holds for every solution $(q,\qG)$ to the adjoint problem
that can be found as a limit of pairs $(\qt,\qtG)$ as specified in the above proof.
\Erem



\vspace{3truemm}


\Begin{thebibliography}{10}

\bibitem{Barbu}
V. Barbu,
``Nonlinear semigroups and differential equations in Banach spaces'',
Noord-hoff, 
Leyden, 
1976.

\bibitem{Barbutrick}
V. Barbu, 
Necessary conditions for nonconvex distributed control problems 
governed by elliptic variational inequalities, 
{\it J. Math. Anal. Appl.} {\bf 80} (1981) \pier{566-597}.

\bibitem{CahH} 
J. W. Cahn and J. E. Hilliard, 
Free energy of a nonuniform system I. Interfacial free energy, 
{\it J. Chem. Phys.\/}
{\bf 2} (1958) 258-267.

\bibitem{CaCo}
L. Calatroni and P. Colli,
Global solution to the Allen--Cahn equation with singular potentials and dynamic boundary conditions,
{\it Nonlinear Anal.\/} {\bf 79} (2013) 12-27.

\pier{ 
\bibitem{CMZ}
 L.\ {C}herfils, A.\ {M}iranville and S.\ {Z}elik,
The {C}ahn--{H}illiard equation with logarithmic potentials,
{\it Milan J. Math.\/} {\bf 79} (2011) 561-596.%
}

\bibitem{CFGS1}
P. Colli, M.H. Farshbaf-Shaker, G. Gilardi and J. Sprekels,
Optimal boundary control of a viscous Cahn--Hilliard system with dynamic boundary condition and double obstacle potentials,
WIAS Preprint No.~2006 (2014) 1-29.

\bibitem{CFS}
P. Colli, M.H. Farshbaf-Shaker and J. Sprekels,
A deep quench approach to the optimal control of an Allen--Cahn equation
with dynamic boundary conditions and double obstacles, 
{\it Appl. Math. Optim.\/} \pier{{\bf 71} (2015) 1-24.}

\pier{
\bibitem{CF} 
P.\ {C}olli and T.\ {F}ukao, 
{C}ahn--{H}illiard equation with dynamic boundary conditions 
and mass constraint on the boundary, 
preprint arXiv:1412.1932 [math.AP] (2014) 1-26.%
}

\bibitem{CGS}
P. Colli, G. Gilardi and J. Sprekels,
On the Cahn--Hilliard equation with dynamic 
boundary conditions and a dominating boundary potential,
{\it J. Math. Anal. Appl.\/} {\bf 419} (2014) 972-994.

\bibitem{CGSvisc}
P. Colli, G. Gilardi and J. Sprekels,
A boundary control problem
for the viscous Cahn-Hilliard equation
with dynamic boundary conditions,
\pier{preprint arXiv:1407.3916 [math.AP] (2014) 1-27.}

\pier{
\bibitem{CS}
P. Colli and J. Sprekels,
Optimal control of an Allen--Cahn equation 
with singular potentials and dynamic boundary condition,
{\it SIAM \juerg{J.} Control Optim.\/} {\bf 53} (2015) 213-234.%
}

\bibitem{EllSt} 
C. M. Elliott and A. M. Stuart, 
Viscous Cahn--Hilliard equation. II. Analysis, 
{\it J. Differential Equations\/} 
{\bf 128} (1996) 387-414.

\bibitem{EllSh} 
C. M. Elliott and S. Zheng, 
On the Cahn--Hilliard equation, 
{\it Arch. Rational Mech. Anal.\/} 
{\bf 96} (1986) 339-357.

\bibitem{GiMiSchi} 
G. Gilardi, A. Miranville and G. Schimperna,
On the Cahn--Hilliard equation with irregular potentials and dynamic boundary conditions,
{\it Commun. Pure Appl. Anal.\/} 
{\bf 8} (2009) 881-912.

\pier{
\bibitem{HW1}
M. Hinterm\"uller and D. Wegner, Distributed optimal control of the 
Cahn--Hilliard system including the case of a double-obstacle 
homogeneous free energy density, {\it SIAM J.
Control Optim.} {\bf 50} (2012) 388-418.
\bibitem{HW2}
M. Hinterm\"uller and D. Wegner,	Optimal control of a semi-discrete 
Cahn--Hilliard--Navier--Stokes system, {\it SIAM J. Control Optim.} 
{\bf 52} (2014) 747-772.
\bibitem{Lions}
J.-L.~Lions,
``Quelques m\'ethodes de r\'esolution des probl\`emes
aux limites non lin\'eaires'',
Dunod; Gauthier-Villars, Paris, 1969.
\bibitem{PRZ} 
J. Pr\"uss, R. Racke and S. Zheng, 
Maximal regularity and asymptotic behavior of solutions for the Cahn--Hilliard equation with dynamic boundary conditions,  
{\it Ann. Mat. Pura Appl.~(4)\/}
{\bf 185} (2006) 627-648.
\bibitem{RZ} 
R. Racke and S. Zheng, 
The Cahn--Hilliard equation with dynamic boundary conditions, 
{\it Adv. Differential Equations\/} 
{\bf 8} (2003) 83-110.
\bibitem{RoSp}
E. Rocca and J. Sprekels,
Optimal distributed control of a nonlocal convective Cahn--Hilliard equation 
by the velocity in 3D, WIAS Preprint No. 1942 (2014), pp. 1-24. 
\bibitem{Simon}
J. Simon,
{Compact sets in the space $L^p(0,T; B)$},
{\it Ann. Mat. Pura Appl.~(4)\/} 
{\bf 146} (1987) 65-96.
\bibitem{Wang}
Q.-F. Wang, Optimal distributed control of nonlinear {C}ahn-{H}illiard 
systems with computational realization, 
{\it J. Math. Sci. (N. Y.)\/} {\bf 177} (2011) 440-458.
\bibitem{WaNa}
Q.-F. Wang and S.-I. Nakagiri, Optimal control of distributed parameter 
system given by Cahn--Hilliard equation, {\it Nonlinear Funct. Anal. Appl.} {\bf 19} (2014) 19-33.
\bibitem{WZ} H. Wu and S. Zheng,
Convergence to equilibrium for the Cahn--Hilliard equation with dynamic boundary conditions, {\it J. Differential Equations\/}
{\bf 204} (2004) 511-531.
\bibitem{ZL1}
X. P. Zhao and C. C. Liu, Optimal control of the convective Cahn--Hilliard equation, 
{\it Appl. Anal.\/} {\bf 92} (2013) 1028-1045.
\bibitem{ZL2}
X. P. Zhao and C. C. Liu, Optimal control for the convective Cahn--Hilliard equation 
in 2D case, {\it Appl. Math. Optim.} {\bf 70} (2014) 61-82.
\bibitem{ZW}
J. Zheng and Y. Wang, Optimal control problem for {C}ahn--{H}illiard equations
with state constraint, {\it J. Dyn. Control Syst.\/}, {\bf 21} (2015) 257-272.%
}
\End{thebibliography}

\End{document}


\noindent
{\bf In the source file after {\tt\expandafter\string\ }\unskip end$\{$document$\}$, 
one can find more references taken from~\cite{CGSvisc}.}

\bibitem{Brezis}
H. Brezis,
``Op\'erateurs maximaux monotones et semi-groupes de contractions
dans les espaces de Hilbert'',
North-Holland Math. Stud.
{\bf 5},
North-Holland,
Amsterdam,
1973.

\bibitem{BreGil}
F. Brezzi and G. Gilardi,
Part~1:
Chapt.~2, Functional spaces, 
Chapt.~3, Partial differential equations,
in ``Finite element handbook'',
H. Kardestuncer and D. H.\ Norrie  eds.,
McGraw-Hill Book Company,
NewYork,
1987.

\bibitem{CFP} 
R. Chill, E. Fa\v sangov\'a and J. Pr\"uss,
Convergence to steady states of solutions of the Cahn--Hilliard equation with dynamic boundary conditions,
{\it Math. Nachr.\/} 
{\bf 279} (2006) 1448-1462.

\bibitem{CF} 
P. Colli and T. Fukao,
The Allen-Cahn equation with dynamic boundary conditions and mass constraints,
preprint arXiv:1405.0116~[math.AP] (2014), pp.~1-23.

\bibitem{CGPS}
P. Colli, G. Gilardi, P. Podio-Guidugli and J. Sprekels, 
Distributed optimal control of a nonstandard system of phase field equations, {\it Cont. Mech. Thermodyn.\/}
{\bf 24} (2012) 437-459.

\bibitem{CGScon}
P. Colli, G. Gilardi and J. Sprekels,
Analysis and boundary control of a nonstandard system of phase field equations, 
{\it Milan J. Math.\/} {\bf 80} (2012) 119-149.

\bibitem{Is} 
H.\ Israel,
Long time behavior of an {A}llen-{C}ahn type equation with a
singular potential and dynamic boundary conditions,
{\it J. Appl. Anal. Comput.} {\bf 2} (2012) 29-56.

\bibitem{LioMag}
J.-L. Lions and E. Magenes,
``Non-homogeneous boundary value problems and applications'',
Vol.~I,
Springer, Berlin, 1972.

\bibitem{MirZelik} 
A. Miranville and S. Zelik,
Robust exponential attractors for Cahn--Hilliard type equations with singular potentials,
{\it Math. Methods Appl. Sci.\/} 
{\bf 27} (2004) 545--582.

\bibitem{Podio}
P. Podio-Guidugli, 
Models of phase segregation and diffusion of atomic species on a lattice,
{\it Ric. Mat.} {\bf 55} (2006) 105-118.


\bye